\documentclass{amsart}
\usepackage{graphicx}
\newtheorem{theorem}{Theorem}[section]
\newtheorem{lemma}[theorem]{Lemma}
\newtheorem{corollary}[theorem]{Corollary}

\theoremstyle{definition}

\newtheorem{example}[theorem]{Example}

\newtheorem{claim}[theorem]{Claim}

\theoremstyle{remark}
\newtheorem{remark}[theorem]{Remark}

\numberwithin{equation}{section}

\begin{document}

\title{Essential state surfaces for knots and links}

\author{Makoto Ozawa}
\address{Department of Natural Sciences, Faculty of Arts and Sciences, Komazawa University, 1-23-1 Komazawa, Setagaya-ku, Tokyo, 154-8525, Japan}
\curraddr{Department of Mathematics and Statistics, The University of Melbourne
Parkville, Victoria 3010, Australia (until March 2011)}
\email{w3c@komazawa-u.ac.jp, ozawam@unimelb.edu.au (until March 2011)}





\begin{abstract}
We study a canonical spanning surface obtained from a knot or link diagram depending on a given Kauffman state, and give a sufficient condition for the surface to be essential.
By using the essential surface, we can see the triviality and splittability of a knot or link from its diagrams.
This has been done on the extended knot or link class which includes all of semiadequate, homogeneous, and most of algebraic knots and links.
In the process of the proof of main theorem, Gabai's Murasugi sum theorem is extended to the case of nonorientable spanning surfaces.
\end{abstract}

\maketitle

\section{Introduction}\label{s:1}
In 1930, Frankl--Pontrjagin \cite{FP} proved the existence of a Seifert surface for any knot, and in 1934, Seifert \cite{S} gave an algorithm to construct a Seifert surface from a knot diagram.
Following Seifert's algorithm, we can construct a spanning surface from a knot diagram depending on a given Kauffman state \cite{K}.
In this paper, we give a sufficient condition for the spanning surface to be essential, and by using the essential surface, we show that a knot or link is trivial (resp. split) if and only if the diagram is trivial (resp. split) under the sufficient condition.

Throughout this paper we work in the piecewise linear category.
For the terminology of knot theory, graph theory and 3-manifold theory, we refer to \cite{C}, \cite{D} and \cite{Ma} respectively.
In Section 2, we give the definition and examples of state surfaces and state main results.
We prepare some lemmas in Section 3, one of which extends Gabai's Murasugi sum theorem \cite{G}, and prove main theorems in Section 4.
In Section 5, we list problems which we should study more.
Finally in Section 6, we summarize some recent progress after this paper.

\section{Definitions, Examples and Results}

\begin{figure}[htbp]
	\begin{center}
		\includegraphics[trim=0mm 0mm 0mm 0mm, width=.6\linewidth]{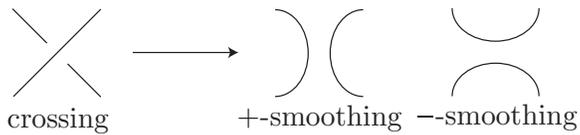}
	\end{center}
	\caption{Two smoothings of a crossing}
	\label{resolution}
\end{figure}

Let $K$ be a knot or link in the 3-sphere $S^3$ and $D$ a connected diagram of $K$ on the 2-sphere $S^2$ which
 separates $S^3$ into two 3-balls, say $B_+, B_-$.
Let $\mathcal{C}=\{c_1,\ldots,c_n\}$ be the set of crossings of $D$.
A map $\sigma : \mathcal{C}\to \{+,-\}$ is called a {\em state} for $D$.
For each crossing $c_i\in \mathcal{C}$, we take a $+$-smoothing or $-$-smoothing according to $\sigma(c_i)=+$ or $-$.
See Figure \ref{resolution}.\footnote{Historically, $+$-smoothing and $-$-smoothing are called $A$-splice and $B$-splice in the most of all papers. It seems to be reasonable that we call these smoothings $+$-smoothing and $-$-smoothing since if we orient a crossing locally so that it has a $\pm$-sign, then a smoothing along the orientation coincides with $\pm$-smoothing.}
Then, we have a collection of loops $l_1,\ldots, l_m$ on $S^2$ and call those {\em state loops}.
Let $\mathcal{L}_{\sigma} = \{l_1,\ldots, l_m\}$ be the set of state loops.

\begin{figure}[htbp]
	\begin{center}
		\includegraphics[trim=0mm 0mm 0mm 0mm, width=.6\linewidth]{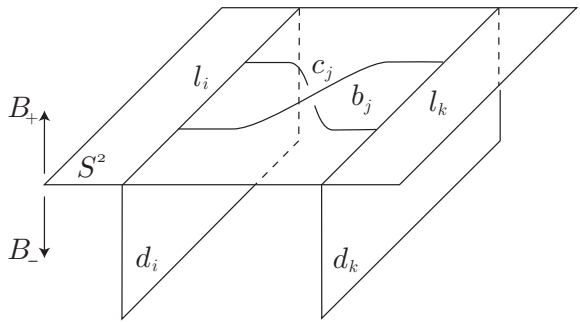}
	\end{center}
	\caption{Recovering a crossing by a half twisted band}
	\label{recover}
\end{figure}

Each state loop $l_i$ bounds a unique disk $d_i$ in $B_-$, and we may assume that these disks are mutually disjoint.
For each crossing $c_j$ and state loops $l_i,l_k$ whose subarcs replaced $c_j$ by $\sigma(c_j)$-smoothing, we attach a half twisted band $b_j$ to $d_i,d_k$ so that it recovers $c_j$.
See Figure \ref{recover} for $\sigma(c_j)=+$.
In this way, we obtain a spanning surface which consists of disks $d_1,\ldots,d_m$ and half twisted bands $b_1,\ldots,b_n$ and call this a {\em $\sigma$-state surface}
 $F_{\sigma}$.

\begin{remark}
Here we mention some historical remarks, which were suggested by J. H. Przytycki.
\begin{enumerate}
\item The state surfaces corresponding to the {\em positive state} $\sigma_+$ (that is, $\sigma_+(c_j)=+$ for all $j$) and {\em negative state} $\sigma_-$ (that is, $\sigma_-(c_j)=-$ for all $j$) were considered for alternating links already in XIX century by Tait (and are called Tait surfaces and nowadays checkerboard surfaces).
\item The state surface corresponding to the {\em Seifert state} $\vec{\sigma}$ (that is, a state determined by an orientation of the knot), which gives the Seifert surface, was introduced by H. Seifert in \cite{S}.
\item Independently, J. H. Przytycki had already thought about the concept of using a surface for any Kauffman state.
See Footnote 2 in \cite{PPS}.
\end{enumerate}

\end{remark}

We may assume that $F_{\sigma}$ intersects $N(K)$ in its collar $N(\partial F_{\sigma};F_{\sigma})$, and overuse the symbol $F_{\sigma}$ instead of $F_{\sigma}\cap E(K)$, where $N(K)$ denotes the regular neighbourhood of $K$ in $S^3$ and $E(K)$ denotes the exterior of $K$.
We take a (twisted) $I$-bundle $F_{\sigma}\tilde{\times} I$ over $F_{\sigma}$ in $E(K)$, and call the associated $\partial I$-bundle $F_{\sigma}\tilde{\times}\partial I$ over $F_{\sigma}$ the {\em interpolating surface} obtained from $F_{\sigma}$ and denote it by $\widetilde{F_{\sigma}}$ since it is a double cover of $F_{\sigma}$.
Note that any interpolating surface $\widetilde{F_{\sigma}}$ is orientable, and it is connected if and only if $F_{\sigma}$ is nonorientable.


We construct a graph $G_{\sigma}$ with signs on edges from $F_{\sigma}$ by regarding a disk $d_i$ as a vertex $v_i$ and a band $b_j$ as an edge $e_j$ which has the same sign $\sigma(c_j)$.
We call the graph $G_{\sigma}$ a {\em $\sigma$-state graph}.
In general, a graph is called a {\em block} if it is connected and has no cut vertex.
It is known that any graph has a unique decomposition into maximal blocks.
Following \cite{LT} and \cite{C1}, we say that a diagram $D$ is {\em $\sigma$-adequate} if $G_{\sigma}$ has no loop,
 and that $D$ is {\em $\sigma$-homogeneous} if in each block of $G_{\sigma}$, all edges have the same sign.
We remark that any diagram of any link is $\sigma$-adequate for some state $\sigma$ (for example, the Seifert state), and $\sigma'$-homogeneous for some state $\sigma'$ (for example, the positive state),
where these states $\sigma$, $\sigma'$ do not coincide generally.

\begin{remark}
As pointed out in \cite{E}, the definition of adequate seems to be slightly different.
See Example \ref{adequate} for the consistency with the original definition.
\end{remark}



\begin{figure}[htbp]
	\begin{center}
		\includegraphics[trim=0mm 0mm 0mm 0mm, width=.9\linewidth]{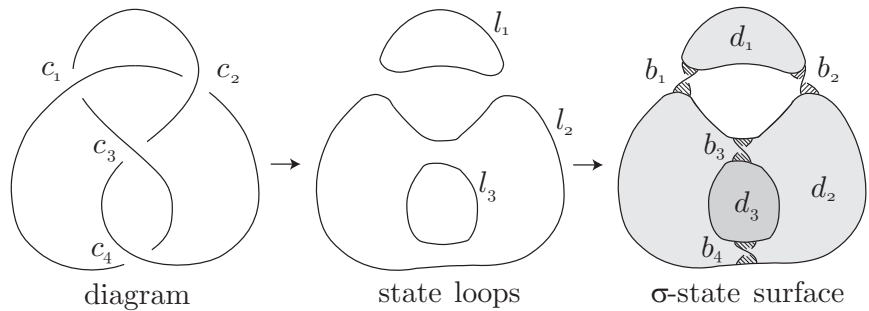}
	\end{center}
	\caption{An example of making a $\sigma$-state surface}
	\label{trefoil}
\end{figure}
\begin{figure}[htbp]
	\begin{center}
		\includegraphics[trim=0mm 0mm 0mm 0mm, width=.4\linewidth]{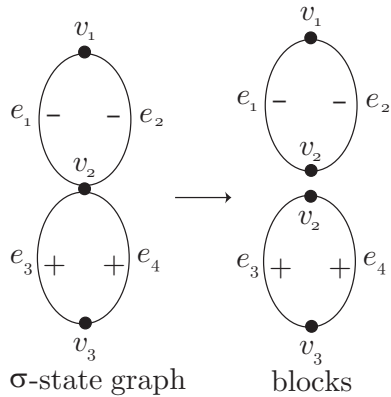}
	\end{center}
	\caption{The corresponding $\sigma$-state graph and its block decomposition}
	\label{figure-8}
\end{figure}
\begin{figure}[htbp]
	\begin{center}
		\includegraphics[trim=0mm 0mm 0mm 0mm, width=.6\linewidth]{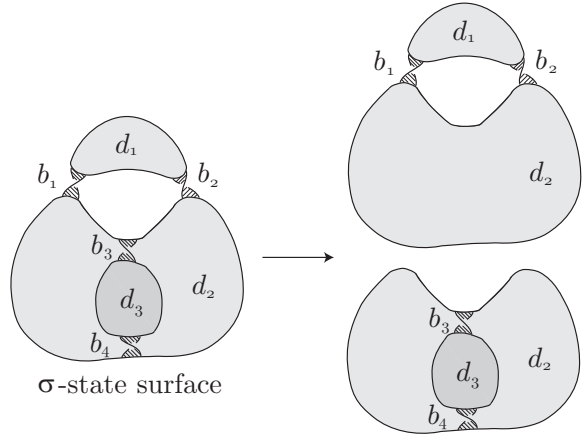}
	\end{center}
	\caption{The corresponding Murasugi decomposition}
	\label{decomposition_fig}
\end{figure}

\begin{example}
Let $D$ be a diagram of the figure eight knot which has $4$ crossings $c_1,c_2,c_3,c_4$ as in Figure \ref{trefoil}.
To make a $\sigma$-state surface, let $\sigma(c_1)=\sigma(c_2)=-$ and $\sigma(c_3)=\sigma(c_4)=+$ for example.
Since the $\sigma$-state graph $G_{\sigma}$ has no loop and all edges in each block have a same sign as in Figure \ref{figure-8}, $D$ is $\sigma$-adequate and $\sigma$-homogeneous.
Moreover, the block decomposition of $G_{\sigma}$ corresponds to a Murasugi decomposition of $F_{\sigma}$. See Figure \ref{decomposition_fig}.
\end{example}

\begin{example}
A diagram $D$ with an orientation is said to be {\em positive} if all crossings have a positive sign.
For any positive diagram $D$, there exists a state $\sigma$ such that $D$ is $\sigma$-adequate and $\sigma$-homogeneous.
Indeed, we can take $\sigma$ so that $\sigma(c_j)=+$ for all $c_j$, namely, the positive state $\sigma_+$.
Also we can take $\sigma$ so that it yields a canonical Seifert surface $F_{\sigma}$, namely, the Seifert state $\vec{\sigma}$.
Note that these states $\sigma_+$ and $\vec{\sigma}$ coincide only on a positive diagram.
\end{example}

\begin{example}
For any alternating diagram $D$ without nugatory crossings, there exist two states $\sigma_1,\sigma_2$ such that $D$ is $\sigma_i$-adequate and $\sigma_i$-homogeneous for $i=1,2$.
Indeed, we can take $\sigma_1=\sigma_+$ (or $\sigma_1=\sigma_-$) and $\sigma_2=\vec{\sigma}$.
\end{example}


\begin{example}
We say that a diagram $D$ is {\em homogeneous} \cite{C1} if $D$ is $\vec{\sigma}$-homogeneous for the Seifert state $\vec{\sigma}$.
Note that $D$ is automatically $\vec{\sigma}$-adequate since the $\vec{\sigma}$-state surface $F_{\vec{\sigma}}$ is orientable and thus $G_{\vec{\sigma}}$ has no loop.

\label{adequate}
We say that a diagram $D$ is {\em semiadequate} \cite{LT} if $D$ is $\sigma$-adequate for the positive state $\sigma_+$ or the negative state $\sigma_-$.
Note that $D$ is automatically $\sigma_{\pm}$-homogeneous since $\sigma_{\pm}(c_j)=\pm$ for all $j$.

We say that a diagram $D$ is {\em adequate} \cite{T} if $D$ is $\sigma$-adequate for both of the positive state $\sigma_+$ and the negative state $\sigma_-$.
Note also that $D$ is automatically $\sigma_{\pm}$-homogeneous since $\sigma_{\pm}(c_j)=\pm$ for all $j$.
\end{example}

\begin{example}\label{arborescent}
We say that an arborescent link $L$ is {\em strictly arborescent} if the absolute value of each weight is greater than $1$.
Note that there exists a diagram $D$ of $L$ and a state $\sigma$ such that $D$ is $\sigma$-adequate and $\sigma$-homogeneous.
Indeed, a strictly arborescent link $L$ is the boundary of a $\sigma$-state surface which is a Murasugi sum of twisted annuli or M\"{o}bius bands with one or more full twists.
See \cite{Ga} or \cite{BS} for the definition and the construction of surfaces for arborescent links.
\end{example}

We review the definition of essential surfaces.

Let $M$ be an orientable compact 3-manifold, $F$ a compact surface properly embedded in $M$, possibly with boundary, except for a 2-sphere, and let $i$ denote the inclusion map $F\to M$.
We say that $F$ is {\em $\pi_1$-injective} if the induced map $i_*:\pi_1(F)\to \pi_1(M)$ is injective, and that $F$ is {\em $\partial$-$\pi_1$-injective} if the induced map $i_*:\pi_1(F,\partial F)\to \pi_1(M,\partial M)$ is injective for every choice of two base points in $\partial F$.
A surface $F$ in $M$ is {\em $\pi_1$-essential} if $F$ is $\pi_1$-injective, $\partial$-$\pi_1$-injective and not $\partial$-parallel in $M$.


A disk $D$ embedded in $M$ is a {\em compressing disk} for $F$ if $D\cap F=\partial D$ and $\partial D$ is an essential loop in $F$.
A disk $D$ embedded in $M$ is {\em $\partial$-compressing disk} for $F$ if $D\cap F\subset \partial D$ is an essential arc in $F$ and $D\cap \partial M=\partial D-\rm{int}(D\cap F)$.
We say that $F$ is {\em incompressible} (resp. {\em $\partial$-incompressible}) if there exists no compressing disk (resp. $\partial$-compressing disk) for $F$.
A surface $F$ in $M$ is {\em essential} if $F$ is incompressible, $\partial$-incompressible and not $\partial$-parallel in $M$.

We remark that a $\sigma$-state surface $F_{\sigma}$ is $\pi_1$-essential in $E(K)$ if and only if the interpolating surface $\widetilde{F_{\sigma}}$ obtained from $F_{\sigma}$ is essential in $E(K)$.

The following main theorem gives a sufficient condition for the state surface to be $\pi_1$-essential.

\begin{theorem}\label{main}
If a diagram is both $\sigma$-adequate and $\sigma$-homogeneous for some state $\sigma$, then the $\sigma$-state surface is $\pi_1$-essential.
\end{theorem}

If $F_{\sigma}$ is nonorientable and $\pi_1$-essential, then the interpolating surface $\widetilde{F_{\sigma}}$ is connected and essential.
Therefore, the knot satisfies the Neuwirth conjecture \cite{N}, which states that for any nontrivial knot $K$, there exists a closed surface $S$ containing $K$ such that $S\cap E(K)$ is connected and essential in $E(K)$.

\begin{corollary}\label{Neuwirth}
If a diagram is $\sigma$-adequate and $\sigma$-homogeneous for a state $\sigma$ except the Seifert state $\vec{\sigma}$, then the knot satisfies the Neuwirth conjecture.
In particular, adequate knots satisfy the Neuwirth conjecture.
\end{corollary}

\begin{remark}
It can be confirmed that every 10 crossing knot diagram in the Rolfsen knot table \cite{R} except for $8_{19}$, $10_{124}$, $10_{128}$, $10_{134}$, $10_{139}$ and $10_{142}$ is $\sigma$-adequate and $\sigma$-homogeneous for a positive or negative state $\sigma$ distinct from the Seifert state $\vec{\sigma}$, and that every 11 crossing knot diagram in the Hoste-Thistlethwaite knot table \cite{HT2} except for $K11_n93$, $K11_n95$, $K11_n118$, $K11_n126$, $K11_n136$, $K11_n169$, $K11_n171$, $K11_n180$ and $K11_n181$ is also $\sigma$-adequate and $\sigma$-homogeneous for a positive or negative state $\sigma$ distinct from the Seifert state $\vec{\sigma}$.
Furthermore, it can be checked that $10_{134}$, $10_{142}$, $K11_n93$, $K11_n95$, $K11_n136$, $K11_n169$, $K11_n171$, $K11_n180$ and $K11_n181$ bound $\pi_1$-essential nonorientable checkerboard surfaces.
(You might need to deform the diagram by the Reidemeister move of type III.)
\end{remark}

\begin{remark}
Futer--Kalfagianni--Purcell \cite{FKP} gave an estimate for the hyperbolic volume of adequate knots by using the guts of state surfaces.
\end{remark}

\begin{remark}
We can construct a spanning surface other than $F_{\sigma}$ from a state $\sigma$ by letting the loop $l_i$ bound a disk $d_i'$ in $B_+$.
Theorem \ref{main} holds for all state surfaces obtained by such a method.
Moreover we can construct a branched surface as in \cite{HT} which consists of disks $d_1,\ldots,d_m$ in $B_-$ and disks $d_1',\ldots,d_m'$ in $B_+$ bounded by $l_1,\ldots,l_m$ respectively, and half twisted bands $b_1,\ldots,b_n$.
\end{remark}

\begin{remark}
Suppose that a diagram $D$ is $\sigma$-adequate and $\sigma$-homogeneous for a state $\sigma$.
If $F_{\sigma}$ is orientable, then it is a minimal genus Seifert surface
 by \cite[Theorem 2]{G} or \cite[Corollary 4.1]{C1}.
On the other hand, M. Hirasawa pointed out that a similar phenomenon need not occur in the nonorientable case. Indeed, there exist 2-bridge links with two continued fractions $-3$ $2$ $-2$ $3$ and $2$ $3$ $2$, where the notation follows Adams's knot book \cite{Adams}.
\end{remark}

\begin{remark}
The converse of Theorem \ref{main} does not hold generally.
It is true that if a $\sigma$-state surface $F_{\sigma}$ is $\pi_1$-essential, then the diagram $D$ is $\sigma$-adequate.
However, in general, it is not true that if a $\sigma$-state surface $F_{\sigma}$ is $\pi_1$-essential, then the diagram $D$ is $\sigma$-homogeneous.
\end{remark}

By using a $\pi_1$-essential state surface, we can show the next theorem which assure us that we can see the triviality and splittability of a knot or link from its diagram.
In this paper, we say that a diagram $D$ is {\em nontrivial} if it contains at least one crossing, and that $D$ is {\em nonsplit} if is is connected.

\begin{theorem}\label{nontrivial}
Let $K$ be a knot or link which admits a $\sigma$-adequate and $\sigma$-homogeneous diagram $D$ without nugatory crossings for some state $\sigma$.
Then,
\begin{enumerate}
	\item[(1)] $D$ is nontrivial if and only if $K$ is nontrivial.
	\item[(2)] $D$ is nonsplit if and only if $K$ is nonsplit.
\end{enumerate}
\end{theorem}

The determining problem for the triviality and splittability was solved about the following classes.
For the triviality, alternating knots \cite{MT}, 
homogeneous links \cite{C1}, semiadequate links \cite{T} and Montesinos knots \cite{LT}.
For the splittability, alternating links \cite{M}, 
homogeneous links \cite{C1}, semiadequate links \cite{T} and positive links \cite{MO2}.

\begin{figure}[htbp]
	\begin{center}
		\includegraphics[trim=0mm 0mm 0mm 0mm, width=.7\linewidth]{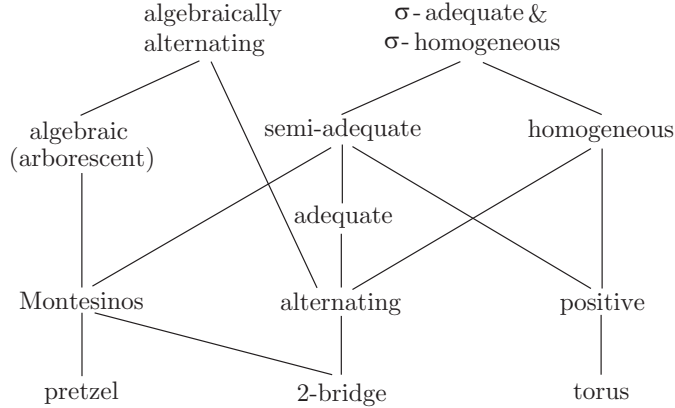}
	\end{center}
	\caption{The Hasse diagram for the set of knot diagrams partially ordered by inclusion}
	\label{hasse}
\end{figure}

\begin{figure}[htbp]
	\begin{center}
		\includegraphics[trim=0mm 0mm 0mm 0mm, width=.5\linewidth]{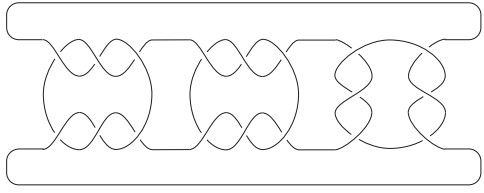}
	\end{center}
	\caption{An example of a diagram of an algebraic link which is not $\sigma$-adequate and $\sigma$-homogeneous for any state $\sigma$}
	\label{algebraic}
\end{figure}

The Hasse diagram of various classes of knots and links are illustrated in Figure \ref{hasse}.
Here, almost all algebraic links have $\sigma$-adequate and $\sigma$-homogeneous diagrams for some state $\sigma$ (see Example \ref{arborescent}), but some algebraic links seem to be not $\sigma$-adequate and $\sigma$-homogeneous for any state $\sigma$ (see Figure \ref{algebraic}).
Algebraically alternating knots and links are defined in \cite{MO3} so that they include both alternating and algebraic knots and links, and some results on closed incompressible surfaces are obtained.

\section{Lemmas}\label{s:2}
The next lemma is stated for knots.
However, it will also hold for a link $K$, so long as $E(K)-F$ is irreducible.
Note that for a connected diagram $D$ and a state surface $F_\sigma$, $E(K)-F_\sigma$ will be a handlebody, hence irreducible.

\begin{lemma}[{\cite[Lemma 2]{MO}}]\label{annulus}
Let $K$ be a knot in $S^3$ and $F$ an incompressible orientable surface properly embedded in $E(K)$.
If $F$ is $\partial$-compressible in $E(K)$, then $F$ is a $\partial$-parallel annulus.
\end{lemma}


Similarly we have:

\begin{lemma}[{\cite[Lemma 2.2]{OT}}]
Let $K$ be a knot in $S^3$ and $F$ a $\pi_1$-injective nonorientable surface properly embedded in $E(K)$.
If $F$ is not $\partial$-$\pi_1$-injective, then $F$ is an unknotted, half-twisted M\"{o}bius band and $K$ is trivial.
\end{lemma}

\begin{lemma}[{\cite[Theorem 9.8]{A}}, {\cite[Proposition 2.3]{MT2}}, {\cite[Theorem 2, 3]{MO}}]\label{checkerboard}
Let $D$ be a reduced, prime, alternating diagram.
Then the checkerboard surface obtained from $D$ is $\pi_1$-essential.
\end{lemma}



Let $F$ be a spanning surface for a link $K$.
Suppose that there exists a 2-sphere $S$ decomposing $S^3$ into two 3-balls $B_1,B_2$ such that $F\cap S$ is a disk.
Put $F_i=F\cap B_i$ for $i=1,2$.
Then we say that $F$ has a {\em Murasugi decomposition} into $F_1$ and $F_2$ and we denote by $F=F_1*F_2$.
Conversely, we say that $F$ is obtained from $F_1$ and $F_2$ by a {\em Murasugi sum} along a disk $F\cap S$.

Put $E=S-int(F\cap S)$ and let $\delta$ be a disk in $B_1$ such that $\delta\cap (F_1\cup E)=\partial \delta\cap (F_1\cup E)=\partial \delta$ and $\partial \delta\cap E$ consists of mutually disjoint arcs $\alpha_1,\ldots,\alpha_n$.
Then, there exist mutually disjoint arcs $\alpha_1',\ldots,\alpha_n'$ in $F\cap S$ which form mutually disjoint loops $\alpha_1\cup \alpha_1',\ldots,\alpha_n\cup\alpha_n'$ in $S$, and there exist mutually disjoint disks $\delta_1',\ldots,\delta_n'$ in $B_2$ which are bounded by $\alpha_1\cup \alpha_1',\ldots,\alpha_n\cup\alpha_n'$ respectively.
We call a disk $\delta\cup (\delta_1'\cup\ldots\cup\delta_n')$ the {\em extended disk} of $\delta$ toward $B_2$.
We remark that the extended disk of $\delta$ is uniquely determined by $\delta$ and generally it intersects $F_2$ in the interior.

The following key lemma extends \cite[Theorem 1]{G} to nonorientable surfaces.

\begin{lemma}\label{decomposition}
If $F_1$ and $F_2$ are $\pi_1$-essential, then $F=F_1*F_2$ is also $\pi_1$-essential.
\end{lemma}

\begin{proof}
We will show that the interpolating surface $\widetilde{F}=F\tilde{\times}\partial I$ is essential.
We note that by \cite[Claim 9]{MO}, $\widetilde{F}$, $\widetilde{F_1}$ and $\widetilde{F_2}$ are incompressible and $\partial$-incompressible in $F\tilde{\times} I$, $F_1\tilde{\times} I$ and $F_2\tilde{\times} I$ respectively.

Let $C$ be a compressing disk for $\widetilde{F}$ in the outside of $F\tilde{\times} I$.
Put $E=S-int(F\cap S)$.
We may assume that $C$ and $E$ are in general position, and that the number of components of $C\cap E$ is minimal over all compressing disks $C$.
If $C\cap E=\emptyset$, then $C$ is a compressing disk for $\widetilde{F_1}$ or $\widetilde{F_2}$.
Otherwise, $C\cap E$ consists of arcs, say $\alpha_1,\ldots,\alpha_p$, and let $\delta_1,\ldots,\delta_q$ be subdisks on $C$ separated by $\alpha_1\cup\cdots\cup\alpha_p$.
For each arc $\alpha_k$, put $\partial \alpha_k=a_k^+\cup a_k^-$.
A subarc $N(a_k^{\pm};\partial C)$ runs over the disk $F\cap S$ and $F-S$.
Then, we mark $a_k^{\pm}$ with an arrow so that it runs from $F\cap S$ to $F-S$.
See Figure \ref{mark}.

\begin{figure}[htbp]
	\begin{center}
		\includegraphics[trim=0mm 0mm 0mm 0mm, width=.9\linewidth]{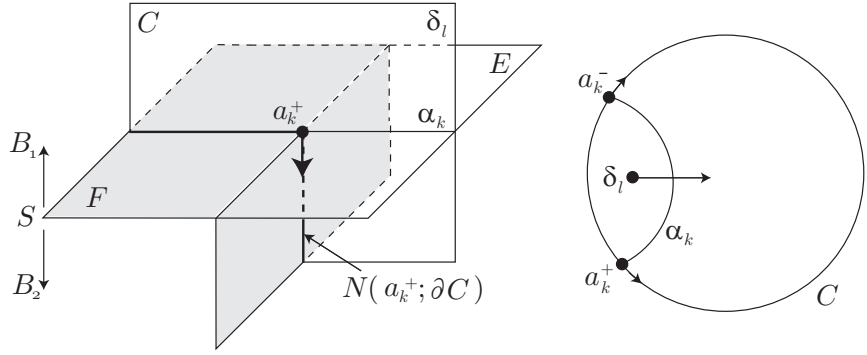}
	\end{center}
	\caption{Marking $a_k^{\pm}$ with an arrow and an induced orientation by $\alpha_k$}
	\label{mark}
\end{figure}

\begin{claim}\label{outermost}
For an outermost arc $\alpha_k$ and the corresponding outermost disk $\delta_l$, both arrows at $a_k^{\pm}$ turn out from $\delta_l$ (as in the right side of Figure \ref{mark}).
\end{claim}

\begin{proof}
Without loss of generality, we may assume that $\delta_l\subset B_1$.
First, suppose that both arrows at $a_k^{\pm}$ turn into $\delta_l$
 (see Figure \ref{arrow1}).
There exists an arc $\alpha_k'$ which connects $a_k^+$ and $a_k^-$ on $F\cap S$, and the loop $\alpha_k\cup \alpha_k'$ bounds a disk $\delta_l'$ in $B_2$.
Then, the extended disk $\delta_l\cup \delta_l'$ toward $B_2$ is a compressing disk for $\widetilde{F_1}$ since we assumed that the number of components of $C\cap E$ is minimal.

Next, suppose that one arrow at $a_k^{\pm}$ turns into $\delta_l$ and another turns out from $\delta_l$
 (see Figure \ref{arrow2}).
Similarly, there exists an arc $\alpha_k'$ which connects $a_k^+$ and $a_k^-$ on $F\cap S$, and the loop $\alpha_k\cup \alpha_k'$ bounds a disk $\delta_l'$ in $B_2$.
Then, the extended disk $\delta_l\cup \delta_l'$ toward $B_2$ is a $\partial$-compressing disk for $\widetilde{F_1}$ since we assumed that the number of components of $C\cap E$ is minimal.
In either case, we have a contradiction.
\end{proof}

\begin{figure}[htbp]
	\begin{center}
		\includegraphics[trim=0mm 0mm 0mm 0mm, width=.8\linewidth]{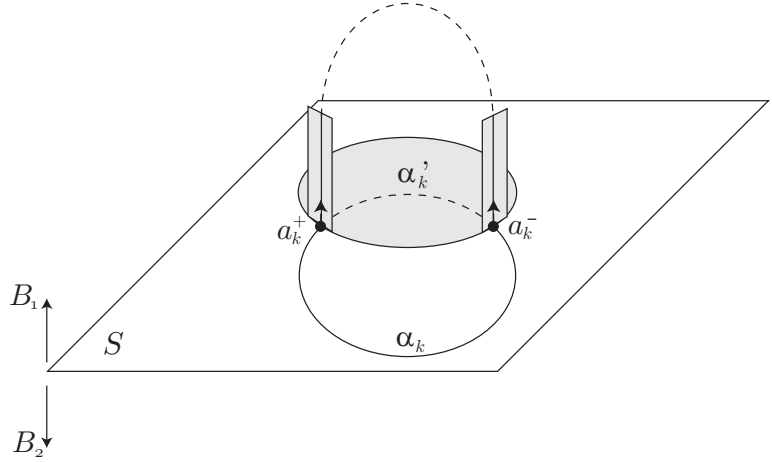}
	\end{center}
	\caption{Both arrows at $a_k^{\pm}$ turn into $\delta_l$}
	\label{arrow1}
\end{figure}

\begin{figure}[htbp]
	\begin{center}
		\includegraphics[trim=0mm 0mm 0mm 0mm, width=.8\linewidth]{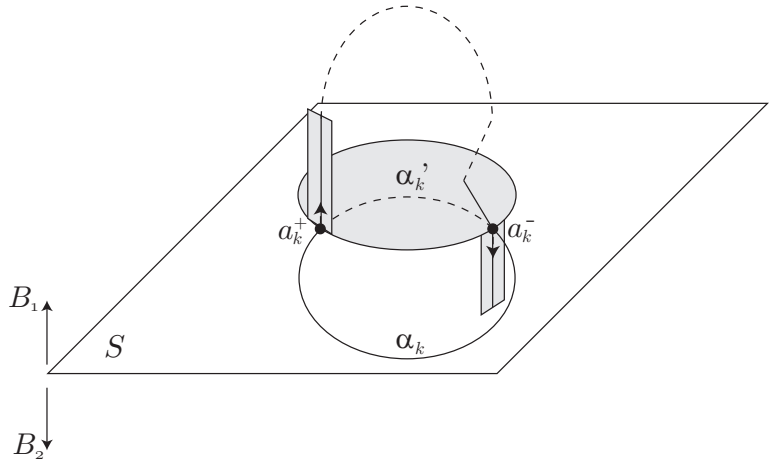}
	\end{center}
	\caption{One arrow at $a_k^+$ turns into $\delta_l$ and another arrow at $a_k^-$ turns out from $\delta_l$}
	\label{arrow2}
\end{figure}

We construct a graph $G$ on $C$ as follows.
We assign a vertex $v_l$ to each subdisk $\delta_l$, and connect two vertices by an edge $e_k$ if the two corresponding subdisks have a common arc $\alpha_k$ of $C\cap E$.
Note that $G$ is a tree, since any arc $\alpha_k$ separates $\delta$.
Since by Claim \ref{outermost}, both arrows at the boundary of an outermost arc are turn out from the corresponding outermost disk, we can assign an orientation to the corresponding outermost edge naturally.
We call such orientation of an edge $e_k$ an {\em induced orientation} by $\alpha_k$.
See Figure \ref{mark}.

A vertex of $G$ has {\em depth $x$} if it becomes a degree $1$ or $0$ vertex after removing all vertices having depth less than $x$, where $x$ is a natural number.
We define vertices corresponding to outermost subdisks as depth $1$.
See Figure \ref{depth}, where the depth of each vertex is indicated.

\begin{figure}[htbp]
	\begin{center}
		\includegraphics[trim=0mm 0mm 0mm 0mm, width=.7\linewidth]{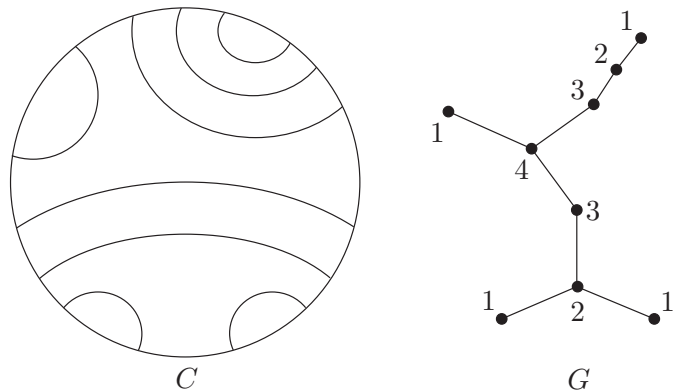}
	\end{center}
	\caption{An example of $C\cap E$ on $C$ and the corresponding graph $G$}
	\label{depth}
\end{figure}

\begin{claim}\label{second}
Every edge of $G$ has an induced orientation and every vertex has an edge oriented outward.
\end{claim}

\begin{proof}
We prove this by induction on the depth of $v_l$.
In case of depth $1$, it was shown by Claim \ref{outermost}.
Next, suppose that Claim \ref{second} holds for vertices having depth less than $x$, and $v_l$ has depth $x$.
Let $N_{< x}(v_l)$ be the set of vertices adjacent to $v_l$ and having depth less than $x$.
Since $G$ has no cycle, any vertex in $N_{<x}(v_l)$ has an edge oriented outward to $v_l$.
Without loss of generality, we may assume that $\delta_l\subset B_1$.
In case $v_l$ becomes a degree $0$ vertex after removing all vertices having depth less than $x$, the extended disk $\delta_l'$ of $\delta_l$ toward $B_2$ is a compressing disk for $\widetilde{F_1}$.
In case $v_l$ becomes a degree $1$ vertex after removing all vertices having depth less than $x$, let $e_k$ be the edge connecting $v_l$ to a vertex except for $N_{<x}(v_l)$, and $\alpha_k$ be the corresponding arc.
First, suppose that both arrows at $a_k^{\pm}$ turn into $\delta_l$.
Then, the extended disk $\delta_l'$ of $\delta_l$ toward $B_2$ is a compressing disk for $\widetilde{F_1}$.
Next, suppose that one arrows at $a_k^{\pm}$ turns into $\delta_l$ and another turns out from $\delta_l$.
Then, the extended disk $\delta_l'$ of $\delta_l$ toward $B_2$ is a $\partial$-compressing disk for $\widetilde{F_1}$.
In either case, we have a contradiction.
Hence, $e_k$ has an induced orientation by $\alpha_k$, and $v_l$ has an edge oriented outward.
\end{proof}




Claim \ref{second} leads us to a contradiction since $G$ is a tree.
Hence $\widetilde{F}$ is incompressible.
By an elementary cut-and-paste argument, this shows that $K$ is nonsplit.
If $\widetilde{F}$ is $\partial$-compressible, then by Lemma \ref{annulus}, it is $\partial$-parallel annulus.
Thus $F$ is a not $\partial$-$\pi_1$-injective M\"{o}bius band and hence one of $F_1$ and $F_2$ is also a not $\partial$-$\pi_1$-injective M\"{o}bius band.
This contradicts that both of $F_1$ and $F_2$ are $\pi_1$-essential.
\end{proof}

\section{Proofs of Theorems}\label{s:3}
\begin{proof} (of Theorem \ref{main})
Suppose that a diagram $D$ is $\sigma$-adequate and $\sigma$-homogeneous for some state $\sigma$.
Then the $\sigma$-state graph $G_{\sigma}$ is decomposed into maximal blocks $G_1,\ldots,G_n$ each of which has no loop and all edges in each block have the same sign.
Let $F_1,\ldots, F_n$ be the corresponding $\sigma$-state surfaces with $G_1,\ldots,G_n$.
Then for each $i$, the boundary $\partial F_i$ represents an alternating diagram which is reduced and prime since $G_i$ has no loop and the block decomposition is maximal.
By Lemma \ref{checkerboard}, $F_i$ is $\pi_1$-essential for each $i$, and by Lemma \ref{decomposition}, $F$ is also $\pi_1$-essential.
\end{proof}


\begin{proof} (of Theorem \ref{nontrivial})
Let $K$ be a knot or link which admits a $\sigma$-adequate and $\sigma$-homogeneous diagram $D$ without nugatory crossings.
By Theorem \ref{main}, a $\sigma$-state surface $F_{\sigma}$ is $\pi_1$-essential.

(1) Suppose that $K$ is nontrivial.
Then, any diagram of $K$ has at least one crossing.
Hence, $D$ is nontrivial.
Conversely, suppose that $D$ is nontrivial.
Since $D$ has at least one crossing and does not have nugatory crossings, there exists a component of $F_{\sigma}$ which is not a disk.
This shows that $K$ is nontrivial.

(2) Suppose that $K$ is nonsplit.
Then, any diagram of $K$ is connected.
Hence, $D$ is nonsplit.
Conversely, suppose that $D$ is nonsplit.
Since $D$ is connected, $F_{\sigma}$ is also connected.
It follows from a cut and paste argument on a splitting sphere that $K$ is nonsplit.
\end{proof}

\section{Problems}\label{s:4}
Here, we list the problems that we should solve in the future.

\begin{enumerate}
	\item Show that there exists a knot which has no $\sigma$-adequate and $\sigma$-homogeneous diagram.
	Furthermore, characterize the nature of knots and links which have $\sigma$-adequate and $\sigma$-homogeneous diagrams.
	\item Determine primeness, satelliteness, fiberedness, smallness and tangle decomposability from a given $\sigma$-adequate and $\sigma$-homogeneous diagram.
	\item Show that for a given knot, the number of all $\sigma$-adequate and $\sigma$-homogeneous diagrams without nugatory crossings is finite.
	\item Classify all knots and links which have $\sigma$-adequate and $\sigma$-homogeneous diagrams.
\end{enumerate}

The author believes essential state surfaces to be useful for solving these problems.

\section{Addendum}\label{s:5}
After 
the first submission to Journal of the Australian Mathematical Society on 7 May 2009, there has been some progress concerning the present paper.
We summarize those results here.

In \cite[Theorem 3]{FKP2}, Futer--Kalfagianni--Purcell just cited our main Theorem 1.9 and they use it to verify the Garoufalidis conjecture on a relation between the boundary slopes of a knot and its colored Jones polynomials.

In \cite{AK}, Adams--Kindred also introduce ``layered surfaces'' for link diagrams which are same as state surfaces in the present paper, and showed that if $K$ is an alternating knot, then any spanning surface of $K$ has the same slope as one of the basic layered surfaces for $K$.

In \cite{CT}, Curtis--Taylor showed by using the above result that for an alternating knot the minimal integral boundary slope is given by the signature plus twice the minimum degree of the Jones polynomial and the maximal integral boundary slope is given by the signature plus twice the maximum degree of the Jones polynomial.

Concerning Stoimenow's paper \cite{S}, the nontriviality of semiadequate links also follows our Theorem 1.15, \cite[Theorem 1.1]{S} gives a partial answer to our Problem 3, \cite[Theorem 1.2]{S} gives a partial answer to our Problem 1, and \cite[Question 5.3]{S} is a part of our Problem 2.

\bigskip

\noindent{\bf Acknowledgement.}
The author would like to thank Prof. Kouki Taniyama for suggesting a research on adequate knots and links by a geometrical method, and thank Prof. Jozef H. Przytycki for letting me know the positive or negative state and the Seifert state.
Finally I deeply appreciate many referees' and editors' many efforts for reading this paper.

\bibliographystyle{amsplain}

\end{document}